\documentclass[a4paper,11pt]{amsart}

\usepackage{graphicx}
\usepackage{mathptmx}
\usepackage{amsmath}
\usepackage{amssymb}
\usepackage{enumitem}
\usepackage{xcolor}
\usepackage{pgfplots}

\newmuskip\pFqmuskip

\newcommand*\pFq[6][8]{%
  \begingroup 
  \pFqmuskip=#1mu\relax
  \mathcode`=\string"8000
  \begingroup\lccode`\~=`\,
  \lowercase{\endgroup\let~}\pFqcomma
  F^{#2}_{#3}{\left(\genfrac..{0pt}{}{#4}{#5}\bigg|#6\right)}%
  \endgroup
}
\newcommand{\pFqcomma}{\mskip\pFqmuskip}

\newtheorem{theorem}{Theorem}[section]

\newtheorem{corollary}[theorem]{Corollary}

\begin{document}

\title[]{Probabilistic generalization of Spivey-type relation for degenerate Bell polynomials}

\author{Taekyun  Kim}
\address{Department of Mathematics, Kwangwoon University, Seoul 139-701, Republic of Korea}
\email{tkkim@kw.ac.kr}
\author{Dae San  Kim}
\address{Department of Mathematics, Sogang University, Seoul 121-742, Republic of Korea}
\email{dskim@sogang.ac.kr}


\begin{abstract}
Following Spivey's pivotal discovery of a recurrence relation for Bell numbers, significant research has emerged concerning various generalizations of Bell numbers and polynomials. For example, Kim and Kim established a Spivey-type recurrence relation specifically for degenerate Bell and Dowling polynomials. In this paper, we extend this work by deriving a probabilistic generalization of Spivey-type recurrence relations for both degenerate Bell and degenerate $r$-Bell polynomials.
\end{abstract}

\maketitle

\markboth{\centerline{\scriptsize Probabilistic generalization of Spivey-type relation for degenerate Bell polynomials}}
{\centerline{\scriptsize Taekyun Kim and Dae San Kim}}

\section{Introduction}
The Stirling number of the second kind ${n \brace k}$ enumerates the number of partitiions of a set with $n$ objects into $k$ nonempty disjoint subsets. The sum of the Stirling numbers of the second kind ${n \brace k}$, called the Bell number and denoted by $\phi_{n}$, counts the total number of partitions of a set with $n$ objects. That is, $\phi_{n}=\sum_{k=0}^{n}{n \brace k}$. The Bell polynomial $\phi_{n}(x)=\sum_{k=0}^{n}{n \brace k}x^{k}$ is a natural polynomial extension of the Bell number $\phi_{n}$ (see \cite{1,2,4,9,11}). \par
Spivey found the following recurrence relation on the Bell numbers (see \cite{11}):
\begin{equation*}
\phi_{l+n}=\sum_{k=0}^{n}\sum_{m=0}^{l}\binom{l}{m}{n \brace k}k^{l-m}\phi_{m}.
\end{equation*}
His relation was extended by Gould-Quaintance to the Bell polynomials (see \cite{5}):
\begin{equation*}
\phi_{l+n}(y)=\sum_{k=0}^{n}\sum_{m=0}^{l}\binom{l}{m}{n \brace k}k^{l-m}y^{k}\phi_{m}(y).
\end{equation*}
Recently, in \cite{8}, Kim-Kim discovered the following relation for the degenerate Bell polynomials (see \eqref{6}, \eqref{7}, \eqref{8}):
\begin{equation}
\phi_{l+n,\lambda}(y)=\sum_{k=0}^{n}\sum_{m=0}^{l}\binom{l}{m}{n \brace k}_{\lambda}(k-n\lambda)_{l-m,\lambda}y^{k}\phi_{m,\lambda}(y). \label{1}
\end{equation} \par
Assume that $Y$ is a random variable whose moment generating function exists in a neighborhood of the origin:
\begin{equation}
E\Big[e^{tY}\Big]=\sum_{n=0}^{\infty}E\big[Y^{n}\big]\frac{t^{n}}{n!},\quad (|t|<r),\label{2}
\end{equation}
exists for some $r>0$, (see \cite{3}). \\
Let $(Y_{j})_{j\ge 1}$ be a sequence of mutually independent copies of the random variable $Y$, and let
\begin{equation}
S_{k}=Y_{1}+Y_{2}+\cdots+Y_{k},\ (k\ge 1),\ \mathrm{with}\quad S_{0}=0. \label{3}
\end{equation} \par
The aim of this paper is to generalize the Spivey-type relation for the  degenerate Bell polynomials in \eqref{1} to that for the probabilistic degenerate Bell polynomials associated with $Y$ (see \eqref{13}, \eqref{14}). Indeed, we obtain (see Theorem 2.1):
\begin{align}
& \phi_{l+n,\lambda}^{Y}(y)=\sum_{k=0}^{n}\sum_{m=0}^{l}\binom{l}{m}\frac{y^{k}}{k!}\sum_{l_{1}+\cdots+l_{k}=n}\binom{n}{l_{1},\dots,l_{k}} \label{4} \\
& \qquad \qquad \times E\bigg[(S_{k}-n\lambda)_{l-m,\lambda}\prod_{i=1}^{k}(Y_{i})_{l_{i},\lambda}\bigg]\phi_{m,\lambda}^{Y}(y), \nonumber
\end{align}
where $l_{1},l_{2},\dots,l_{k}$ are positive integers. \par
Furthermore, the relation in \eqref{4} is extended to the probabilistic degenerate $r$-Bell polynomials associated with $Y$ (see \eqref{22}) as follows (see Theorem 2.3):
\begin{align*}
&\phi_{j+n,\lambda}^{(r,Y)}(y)=\sum_{l=0}^{n}\binom{n}{l}(r)_{n-l,\lambda}\sum_{k=0}^{l}\sum_{m=0}^{j}\frac{y^{k}}{k!}\binom{j}{m}\sum_{l_{1}+\cdots+l_{k}=l}\binom{l}{l_{1},\dots,l_{k}} \\
&\qquad\times E\bigg[\big(S_{k}-n\lambda)_{j-m,\lambda}\prod_{i=1}^{k}(Y_{i})_{l_{i},\lambda}\bigg]\phi_{m,\lambda}^{(r,Y)}(y).
\end{align*} \par
The degenerate exponentials are defined as
\begin{equation}
e_{\lambda}^{x}(t)=\sum_{n=0}^{\infty}(x)_{n,\lambda}\frac{t^{n}}{n!}=(1+\lambda t)^{\frac{x}{\lambda}},\quad e_{\lambda}(t)=e_{\lambda}^{1}(t),\quad (\mathrm{see}\ [8,9,12]),\label{5}
\end{equation}
where
\begin{equation*}
(x)_{0,\lambda}=1,\quad (x)_{n,\lambda}=x(x-\lambda)(x-2\lambda)\cdots\big(x-(n-1)\lambda\big),\ (n\ge 1).
\end{equation*}
The degenerate Stirling numbers of the second kind are defined by
\begin{equation}
(x)_{n,\lambda}=\sum_{k=0}^{n}{n \brace k}_{\lambda}(x)_{k},\quad (n\ge 0),\quad (\mathrm{see} \ [9]),\label{6}
\end{equation}
where
\begin{displaymath}
(x)_{0}=1,\quad (x)_{n}=x(x-1)\cdots(x-n+1),\ (n\ge 1).
\end{displaymath}
The degenerate Bell polynomials are given by
\begin{equation}
\phi_{n,\lambda}(x)=\sum_{k=0}^{n}{n \brace k}_{\lambda}x^{k},\quad (n\ge 0),\quad (\mathrm{see}\ [9]).\label{7}
\end{equation}
When $x=1,\ \phi_{n,\lambda}=\phi_{n,\lambda}(1)$ are called the degenerate Bell numbers. Note that
\begin{displaymath}
\lim_{\lambda\rightarrow 0}\phi_{n,\lambda}(x)=\phi_{n}(x),\quad \lim_{\lambda\rightarrow 0}{n \brace k}_{\lambda}={n \brace k}.
\end{displaymath}
From \eqref{7}, we have
\begin{equation}
e^{x(e_{\lambda}(t)-1)}=\sum_{n=0}^{\infty}\phi_{n,\lambda}(x)\frac{t^{n}}{n!},\quad (\mathrm{see}\ [9]).\label{8}
\end{equation} \par
For $r\in\mathbb{N}$ and $n\ge 0$, the $r$-Stirling numbers of the second kind are defined by
\begin{equation}
(x+r)_{n,\lambda}=\sum_{k=0}^{n}{n+r \brace k+r}_{r,\lambda}(x)_{k},\quad (n\ge 0),\quad (\mathrm{see}\ [8]).\label{9}
\end{equation}
In view of \eqref{7}, the degenerate $r$-Bell polynomials are given by
\begin{equation}
\phi_{n,\lambda}^{(r)}(x)=\sum_{k=0}^{n}{n+r \brace k+r}_{r,\lambda}x^{k},\quad (n\ge 0),\quad (\mathrm{see}\ [8]).\label{10}
\end{equation}
When $x=1,\ \phi_{n,\lambda}^{(r)}=\phi_{n,\lambda}^{(r)}(1)$ are called the degenerate $r$-Bell numbers. \par
We note that
\begin{equation}
g(x+t)=\sum_{n=0}^{\infty}\frac{g^{(n)}(x)}{n!}t^{n}=\sum_{n=0}^{\infty}\frac{t^{n}D_{x}^{n}}{n!}g(x)=e^{tD_{x}}g(x),\quad (\mathrm{see}\ [7]), \label{11}
\end{equation}
where $D_{x}=\frac{d}{dx}$. \par
In \cite{9}, the probabilistic degenerate Stirling numbers of the second kind associated with $Y$ are defined by
\begin{equation}
\frac{1}{k!}\Big(E\big[e_{\lambda}^{Y}(t)\big]-1\Big)^{k}=\sum_{n=k}^{\infty}{n \brace k}_{Y,\lambda}\frac{t^{n}}{n!}.\label{12}	
\end{equation}
Thus, by \eqref{12}, we get
\begin{equation*}
{n \brace k}_{Y,\lambda}=\frac{1}{k!}\sum_{l=0}^{k}\binom{k}{l}(-1)^{k-l}E\big[(S_{l})_{n,\lambda}\big],\quad (0\le k\le n).
\end{equation*}
In view of \eqref{8}, the probabilistic degenerate Bell polynomials associated with $Y$ are defined by
\begin{equation}
e^{x\big(E[e_{\lambda}^{Y}(t)]-1\big)}=\sum_{n=0}^{\infty}\phi_{n,\lambda}^{Y}(x)\frac{t^{n}}{n!},\quad (\mathrm{see}\ [9]). \label{13}
\end{equation}
By \eqref{10} and \eqref{12}, we get
\begin{equation}
\phi_{n,\lambda}^{Y}(x)=\sum_{k=0}^{n}{n \brace k}_{Y,\lambda}x^{k}, \quad (n\ge 0),\label{14}
\end{equation}
and
\begin{equation*}
\phi_{n+1,\lambda}^{Y}(x)=x\sum_{k=0}^{n}\binom{n}{k}E\Big[(Y)_{k+1,\lambda}\Big]\phi_{n-k,\lambda}^{Y}(x),\quad (\mathrm{see}\ [9]).
\end{equation*}

\section{Probabilistic generalization of Spivey-type relation for degenerate Bell polynomials}
Let $Y, (Y_{j})_{j\ge 1}$, and $S_{k},\,\,(k \ge 0)$, be as in \eqref{2} and \eqref{3}.
By \eqref{5}, we get
\begin{align}
e_{\lambda}^{Y}(x+t)&=\big(1+\lambda(x+t)\big)^{\frac{Y}{\lambda}}=(1+\lambda x)^{\frac{Y}{\lambda}}\bigg(1+\frac{\lambda t}{1+\lambda x}\bigg)^{\frac{Y}{\lambda}}\label{15}\\
&=e_{\lambda}^{Y}(x)e_{\lambda}^{Y}\bigg(\frac{t}{1+\lambda x}\bigg).\nonumber	
\end{align}
From \eqref{11} and \eqref{15}, we note that
\begin{align}
e^{tD_{x}}e^{y\big(E\big[e_{\lambda}^{Y}(x)\big]-1\big)}&=e^{y\big(E\big[e_{\lambda}^{Y}(x+t)\big]-1\big)}=e^{y\big(E\big[e_{\lambda}^{Y}(x)\big(e_{\lambda}^{Y}\big(\frac{t}{1+\lambda x}\big)-1\big)+e_{\lambda}^{Y}(x)]-1\big)} \label{16}\\
&=e^{yE\big[e_{\lambda}^{Y}(x)\big(e_{\lambda}^{Y}\big(\frac{t}{1+\lambda x}\big)-1\big)\big]}e^{y\big(E[e_{\lambda}^{Y}(x)]-1\big)}. \nonumber
\end{align}

Now, we observe that
\begin{align}
& e^{yE\big[e_{\lambda}^{Y}(x)\big(e_{\lambda}^{Y}\big(\frac{t}{1+\lambda x}\big)-1\big)\big]} =\sum_{k=0}^{\infty}\frac{y^{k}}{k!}\bigg(E\bigg[e_{\lambda}^{Y}(x)\bigg(e_{\lambda}^{Y}\bigg(\frac{t}{1+\lambda x}\bigg)-1\bigg)\bigg]\bigg)^{k}\label{17}\\
&=\sum_{k=0}^{\infty}\frac{y^{k}}{k!}E\bigg[e_{\lambda}^{Y_{1}+\cdots+Y_{k}}(x)\bigg(e_{\lambda}^{Y_{1}}\bigg(\frac{t}{1+\lambda x}\bigg)-1\bigg)\cdots \bigg(e_{\lambda}^{Y_{k}}\bigg(\frac{t}{1+\lambda x}\bigg)-1\bigg)\bigg]\nonumber\\
&=\sum_{k=0}^{\infty}\frac{y^{k}}{k!}\sum_{n=k}^{\infty}\bigg(\sum_{l_{1}+\cdots+l_{k}=n}\binom{n}{l_{1},\dots,l_{k}}E\Big[e_{\lambda}^{S_{k}}(x)(Y_{1})_{l_{1},\lambda}(Y_{2})_{l_{2},\lambda}\cdots (Y_{k})_{l_{k},\lambda}\Big]\bigg)\frac{\big(\frac{t}{1+\lambda x}\big)^{n}}{n!}\nonumber \\
&=\sum_{k=0}^{\infty}\frac{y^{k}}{k!}\sum_{n=k}^{\infty}\bigg(\sum_{l_{1}+\cdots+l_{k}=n}\binom{n}{l_{1},\dots,l_{k}}E\Big[e_{\lambda}^{S_{k}-n\lambda}(x)(Y_{1})_{l_{1},\lambda}(Y_{2})_{l_{2},\lambda}\cdots (Y_{k})_{l_{k},\lambda}\Big]\bigg)\frac{t^{n}}{n!} \nonumber\\
&=\sum_{n=0}^{\infty}\bigg(\sum_{k=0}^{n}\frac{y^{k}}{k!} \sum_{l_{1}+\cdots+l_{k}=n}\binom{n}{l_{1},\dots,l_{k}}E\Big[e_{\lambda}^{S_{k}-n\lambda}(x)\prod_{i=1}^{k}(Y_{i})_{l_{i},\lambda}\Big]\bigg)\frac{t^{n}}{n!},\nonumber
\end{align}
where $l_{1},l_{2},\dots,l_{}$ are positive integers. \par
By \eqref{13}, \eqref{16}, and \eqref{17}, we get
\begin{align}
& e^{tD_{x}} e^{y\big(E[e_{\lambda}^{Y}(x)]-1\big)}= e^{yE\big[e_{\lambda}^{Y}(x)\big(e_{\lambda}^{Y}\big(\frac{t}{1+\lambda x}\big)-1\big)\big]}e^{y\big(E[e_{\lambda}^{Y}(x)] -1\big)}\label{18}\\
&=\sum_{n=0}^{\infty}\sum_{k=0}^{n}\frac{y^{k}}{k!}\sum_{l_{1}+\cdots+l_{k}=n}\binom{n}{l_{1},\dots,l_{k}}\frac{t^{n}}{n!} \nonumber \\
& \quad\quad\quad \times \sum_{j=0}^{\infty}E\bigg[(S_{k}-n\lambda)_{j,\lambda}\prod_{i=1}^{k}(Y_{i})_{l_{i},\lambda}\bigg]\frac{x^{j}}{j!}\sum_{m=0}^{\infty}\phi_{m,\lambda}^{Y}(y)\frac{x^{m}}{m!}\nonumber\\
&=\sum_{n=0}^{\infty}\sum_{k=0}^{n}\frac{y^{k}}{k!}\sum_{l_{1}+\cdots+l_{k}=n}\binom{n}{l_{1},\dots,l_{k}}\frac{t^{n}}{n!} \nonumber \\
& \quad\quad\quad \times \sum_{l=0}^{\infty}\sum_{m=0}^{l}\binom{l}{m}E\bigg[(S_{k}-n\lambda)_{l-m,\lambda}\prod_{i=1}^{k}(Y_{i})_{l_{i},\lambda}\bigg]\phi_{m,\lambda}^{Y}(y)\frac{x^{l}}{l!} \nonumber \\
&=\sum_{n=0}^{\infty}\sum_{l=0}^{\infty}\sum_{k=0}^{n}\sum_{m=0}^{l}\binom{l}{m}\frac{y^{k}}{k!}\sum_{l_{1}+\cdots+l_{k}=n}\binom{n}{l_{1},\dots,l_{k}} \nonumber \\
&\quad\quad\quad \times E\bigg[(S_{k}-n\lambda)_{l-m,\lambda}\prod_{i=1}^{k}(Y_{i})_{l_{i},\lambda}\bigg]\phi_{m,\lambda}^{Y}(y)\frac{t^{n}}{n!}\frac{x^{l}}{l!}. \nonumber
\end{align}
On the other hand, by \eqref{13}, we get
\begin{align}
e^{tD_{x}} e^{y\big(E[e_{\lambda}^{Y}(x)]-1\big)}&=e^{tD_{x}}\sum_{l=0}^{\infty}\phi_{l,\lambda}^{Y}(y)\frac{x^{l}}{l!}\label{19}\\
&=\sum_{n=0}^{\infty}\frac{t^{n}}{n!}D_{x}^{n}\sum_{l=0}^{\infty}\phi_{l,\lambda}^{Y}(y)\frac{x^{l}}{l!}\nonumber\\
&=\sum_{n=0}^{\infty}\sum_{l=0}^{\infty}\phi_{l+n,\lambda}^{Y}(y)\frac{t^{n}}{n!}\frac{x^{l}}{l!}. \nonumber
\end{align}
Therefore, by \eqref{18} and \eqref{19}, we obtain the following theorem.
\begin{theorem}
For $n,l\ge 0$, we have
\begin{equation}
\begin{aligned}
& \phi_{l+n,\lambda}^{Y}(y)=\sum_{k=0}^{n}\sum_{m=0}^{l}\binom{l}{m}\frac{y^{k}}{k!}\sum_{l_{1}+\cdots+l_{k}=n}\binom{n}{l_{1},\dots,l_{k}}\\
& \qquad\times E\bigg[(S_{k}-n\lambda)_{l-m,\lambda}\prod_{i=1}^{k}(Y_{i})_{l_{i},\lambda}\bigg]\phi_{m,\lambda}^{Y}(y),	
\end{aligned} \label{20}
\end{equation}
where $l_{1},l_{2},\dots,l_{k}$ are positive integers.
\end{theorem}
When $Y=1$, from \eqref{20} we have the following relation in \eqref{1}:
\begin{align*}
\phi_{l+n,\lambda}(y)&=\sum_{k=0}^{n}\sum_{m=0}^{l}\binom{l}{m} \big(k-n\lambda\big)_{l-m,\lambda}y^{k}\phi_{m,\lambda}(y)\\
&\quad\quad\quad\quad \times \frac{1}{k!}\sum_{l_{1}+\cdots+l_{k}=n}\binom{n}{l_{1},l_{2},\dots,l_{k}}\prod_{i=1}^{k}(1)_{l_{i},\lambda} \\
&=\sum_{k=0}^{n}\sum_{m=0}^{l}\binom{l}{m}{n \brace k}_{\lambda}(k-n\lambda)_{l-m,\lambda}y^{k}\phi_{m,\lambda}(y).
\end{align*}
Here we note that
\begin{equation*}
{n \brace k}_{\lambda}=\frac{1}{k!}\sum_{l_{1}+\cdots+l_{k}=n}\binom{n}{l_{1},l_{2},\dots,l_{k}}\prod_{i=1}^{k}(1)_{l_{i},\lambda},
\end{equation*}
where $l_{1},l_{2},\dots,l_{k}$ are positive integers. \\
Further, letting $ \lambda \rightarrow 0$ yields Gould-Quaintance extension for Bell polynomials of Spivey's recurrence relation for Bell numbers.
\begin{displaymath}
\phi_{l+n}(y)=\sum_{k=0}^{n}\sum_{m=0}^{l}\binom{l}{m}{n \brace k}k^{l-m}y^{k}\phi_{m}(y),\quad (\mathrm{see}\ [5,7]).
\end{displaymath}

Letting $y=1$ in \eqref{20}, we get the following Spivey-like relation for the probabilistic degenerate Bell numbers associated with $Y$.
\begin{corollary}
For $n,l\ge 0$, we have
\begin{displaymath}
\phi_{l+n,\lambda}^{Y}= \sum_{k=0}^{n}\sum_{m=0}^{l}\binom{l}{m}\frac{1}{k!}\sum_{l_{1}+\cdots+l_{k}=n}\binom{n}{l_{1},\dots,l_{k}} E\bigg[\big(S_{k}-n\lambda)_{l-m,\lambda}\prod_{i=1}^{k}(Y_{i})_{l_{i},\lambda}\bigg]\phi_{m,\lambda}^{Y},
\end{displaymath}
where $l_{1},l_{2},\dots,l_{k}$ are positive integers.
\end{corollary}

For $r\in\mathbb{N}$, we define the probabilistic degenerate $r$-Stirling numbers of the second kind associated with $Y$ as
\begin{equation}
\frac{1}{k!}\Big(E\Big[e_{\lambda}^{Y}(t)\Big]-1\Big)^{k}e_{\lambda}^{r}(t)=\sum_{n=k}^{\infty}{n+r \brace k+r}_{r,\lambda}^{Y}\frac{t^{n}}{n!},\quad (k\ge 0).\label{21}
\end{equation}
When $Y=1$, we have (see \eqref{9})
\begin{displaymath}
\quad {n+r \brace k+r}^{Y}_{r,\lambda}={n+r \brace k+r}_{r,\lambda},\quad (n\ge k\ge 0).
\end{displaymath}
Also, we define the probabilistic degenerate $r$-Bell polynomials associated with $Y$ by
\begin{equation}
e^{x\big(E\big[e_{\lambda}^{Y}(t)\big]-1\big)} e_{\lambda}^{r}(t)=\sum_{n=0}^{\infty}\phi_{n,\lambda}^{(r,Y)}(x)\frac{t^{n}}{n!}. \label{22}
\end{equation}
When $x=1,\ \phi_{n,\lambda}^{(r,Y)}=\phi_{n,\lambda}^{(r,Y)}(1)$ are called the probabilistic degenerate $r$-Bell numbers associated with $Y$. \par
By \eqref{8}, \eqref{15}-\eqref{17}, and \eqref{22}, we get
\begin{align}
&e^{tD_{x}}\Big(e^{y\big(E\big[e_{\lambda}^{Y}(x)\big]-1\big)} e_{\lambda}^{r}(x)\Big)=e^{y\big(E\big[e_{\lambda}^{Y}(x+t)\big]-1\big)}e_{\lambda}^{r}(x+t) \label{23} \\
&= e^{yE\big[e_{\lambda}^{Y}(x)\big(e_{\lambda}^{Y}\big(\frac{t}{1+\lambda x}-1\big)\big]}e_{\lambda}^{r}\bigg(\frac{t}{1+\lambda x}\bigg) e^{y\big(E[e_{\lambda}^{Y}(x)]-1\big)} e_{\lambda}^{r}(x) \nonumber \\
&=\sum_{l=0}^{\infty}\sum_{k=0}^{l}\frac{y^{k}}{k!} \sum_{l_{1}+\cdots+l_{k}=l}\binom{l}{l_{1},\dots,l_{k}}\frac{t^{l}}{l!} E\bigg[e_{\lambda}^{S_{k}-l\lambda}(x)\prod_{i=1}^{k}(Y_{i})_{l_{i},\lambda}\bigg]\nonumber\\
&\quad\times \sum_{s=0}^{\infty}(r)_{s,\lambda}\frac{t^{s}}{s!}e_{\lambda}^{-s \lambda}(x) \sum_{m=0}^{\infty}\phi_{m,\lambda}^{(r,Y)}(y)\frac{x^{m}}{m!}\nonumber \\
&=\sum_{n=0}^{\infty}\sum_{l=0}^{n}\binom{n}{l}(r)_{n-l,\lambda}\sum_{k=0}^{l}\frac{y^{k}}{k!} \sum_{l_{1}+\cdots+l_{k}=l}\binom{l}{l_{1},\dots,l_{k}}\frac{t^{n}}{n!} \nonumber\\
&\quad\times E\bigg[e_{\lambda}^{S_{k}-l\lambda-(n-l)\lambda}(x)\prod_{i=1}^{k}(Y_{i})_{l_{i},\lambda}\bigg]\sum_{m=0}^{\infty}\phi_{m,\lambda}^{(r,Y)}(y)\frac{x^{m}}{m!} \nonumber \\
&=\sum_{n=0}^{\infty}\sum_{l=0}^{n}\binom{n}{l}(r)_{n-l,\lambda}\sum_{k=0}^{l}\frac{y^{k}}{k!} \sum_{l_{1}+\cdots+l_{k}=l}\binom{l}{l_{1},\dots,l_{k}}\frac{t^{n}}{n!} \nonumber\\
&\quad\times \sum_{a=0}^{\infty} E\bigg[(S_{k}-n\lambda)_{a,\lambda}\prod_{i=1}^{k}(Y_{i})_{l_{i},\lambda}\bigg]\frac{x^{a}}{a!}\sum_{m=0}^{\infty}\phi_{m,\lambda}^{(r,Y)}(y)\frac{x^{m}}{m!}\nonumber\\
&=\sum_{n=0}^{\infty}\sum_{l=0}^{n}\binom{n}{l}(r)_{n-l,\lambda}\sum_{k=0}^{l}\frac{y^{k}}{k!} \sum_{l_{1}+\cdots+l_{k}=l}\binom{l}{l_{1},\dots,l_{k}}\frac{t^{n}}{n!} \nonumber\\
&\quad\times \sum_{j=0}^{\infty} \sum_{m=0}^{j}\binom{j}{m}E\bigg[(S_{k}-n\lambda)_{j-m,\lambda}\prod_{i=1}^{k}(Y_{i})_{l_{i},\lambda}\bigg]\phi_{m,\lambda}^{(r,Y)}(y)\frac{x^{j}}{j!}\nonumber\\
&=\sum_{n=0}^{\infty}\sum_{j=0}^{\infty}\sum_{m=0}^{j}\sum_{l=0}^{n}\binom{n}{l}(r)_{n-l,\lambda}\binom{j}{m}\sum_{k=0}^{l}\frac{y^{k}}{k} \sum_{l_{1}+\cdots+l_{k}=l}\binom{l}{l_{1},\dots,l_{k}} \nonumber \\
&\quad \times E\bigg[\Big(S_{k}-n\lambda\Big)_{j-m,\lambda}\prod_{i=1}^{k}(Y_{i})_{l_{i},\lambda}\bigg]\phi_{m,\lambda}^{(r,Y)}(y)\frac{t^{n}}{n!}\frac{x^{j}}{j!}, \nonumber
\end{align}
where $l_{1},l_{2},\dots,l_{k}$ are positive integers. \par
On the other hand, by \eqref{22}, we get
\begin{equation}
\begin{aligned}
&e^{tD_{x}}\Big(e^{y\big(E\big[e_{\lambda}^{Y}(x)\big]-1\big)} e_{\lambda}^{r}(x)\Big)\\
&=\sum_{n=0}^{\infty}\frac{t^{n}}{n!}D_{x}^{n}\sum_{j=0}^{\infty}\phi_{j,\lambda}^{(r,Y)}(y)\frac{x^{j}}{j!}=\sum_{n=0}^{\infty}\sum_{j=0}^{\infty}\phi_{j+n,\lambda}^{(r,Y)}(y)\frac{t^{n}}{n!]}\frac{x^{j}}{j!}.
\end{aligned}\label{24}
\end{equation}
Therefore, by \eqref{23} and \eqref{24}, we obtain the following theorem.
\begin{theorem}
For $n,j\ge 0$, we have
\begin{align*}
&\phi_{j+n,\lambda}^{(r,Y)}(y)=\sum_{l=0}^{n}\binom{n}{l}(r)_{n-l,\lambda}\sum_{k=0}^{l}\sum_{m=0}^{j}\frac{y^{k}}{k!}\binom{j}{m}\sum_{l_{1}+\cdots+l_{k}=l}\binom{l}{l_{1},\dots,l_{k}} \\
&\qquad\times E\bigg[\big(S_{k}-n\lambda)_{j-m,\lambda}\prod_{i=1}^{k}(Y_{i})_{l_{i},\lambda}\bigg]\phi_{m,\lambda}^{(r,Y)}(y),
\end{align*}
where $l_{1},l_{2},\dots,l_{k}$ are positive integers.
\end{theorem}
When $Y=1$, we have
\begin{align}
&\phi_{j+n,\lambda}^{(r)}(y)=\sum_{l=0}^{n}\sum_{m=0}^{j}\binom{n}{l}\binom{j}{m}(r)_{n-l,\lambda}\sum_{k=0}^{l}\frac{y^{k}}{k!}\nonumber\\
&\qquad\times \sum_{l_{1}+\cdots+l_{k}=l}\binom{l}{l_{1},\dots,l_{k}}\prod_{i=1}^{k}(1)_{l_{i},\lambda}(k-n \lambda)_{j-m,\lambda}\phi_{m,\lambda}^{(r)}(y) \label{25}\\
&=\sum_{l=0}^{n}\sum_{m=0}^{j}\binom{n}{l}\binom{j}{m}(r)_{n-l,\lambda}\sum_{k=0}^{l}y^{k}{l \brace k}_{\lambda}(k-n\lambda)_{j-m,\lambda}\phi_{m,\lambda}^{(r)}(y) \nonumber\\
&=\sum_{k=0}^{n}\sum_{m=0}^{j}\sum_{l=k}^{n}\binom{n}{l}(r)_{n-l,\lambda}{l \brace k}_{\lambda}\binom{j}{m}y^{k}(k-n\lambda)_{j-m,\lambda}\phi_{m,\lambda}^{(r)}(y). \nonumber
\end{align}
From \eqref{21}, we note that
\begin{align}
\sum_{n=k}^{\infty}{n+r \brace k+r}_{r,\lambda}\frac{t^{n}}{n!}&=\frac{1}{k!}\big(e_{\lambda}(t)-1\big)^{k}e_{\lambda}^{r}(t) \label{26}\\
&=\sum_{l=k}^{\infty}{l \brace k}_{\lambda}\frac{t^{l}}{l!}\sum_{j=0}^{\infty}(r)_{j,\lambda}\frac{t^{j}}{j!}\nonumber\\
&=\sum_{n=k}^{\infty}\sum_{l=k}^{n}\binom{n}{l}{l \brace k}_{\lambda}(r)_{n-l,\lambda}\frac{t^{n}}{n!}.\nonumber
\end{align}
Thus, by comparing the coefficients on both sides of \eqref{26}, we get
\begin{equation}
{n+ r \brace k+r}_{r,\lambda}=\sum_{l=k}^{n}\binom{n}{l}{l \brace k}_{\lambda}(r)_{n-l,\lambda},\quad (n\ge k\ge 0). \label{27}	
\end{equation}
Therefore, by \eqref{25} and \eqref{27}, we obtain the following corollary.
\begin{corollary}
For $j,n\ge 0$, we have
\begin{equation*}
\phi_{j+n,\lambda}^{(r)}(y)=\sum_{k=0}^{n}\sum_{m=0}^{j}{n+r \brace k+r}_{r,\lambda}\binom{j}{m}\big(k-n\lambda\big)_{j-m,\lambda}y^{k}\phi_{m,\lambda}^{(r)}(y).
\end{equation*}
\end{corollary}

\section{Conclusion}
In 2008, Spivey's discovery of a recurrence relation for the Bell numbers $\phi_{n}$ opened a new avenue of research. Since then, significant work has been dedicated to exploring various Spivey-type recurrence relations (see \cite{5,6,7,8,10}). This includes generalizing Spivey's relation to probabilistic Bell and $r$-Bell polynomials associated with $Y$. Here $Y$ is a random variable whose moment generating function exists in a neighborhood of the origin. \par
In this paper, we further expanded this area by proving Spivey-type recurrence relations for probabilistic degenerate Bell and $r$-Bell polynomials associated with $Y$. These polynomials are degenerate versions of the aforementioned probabilistic Bell and $r$-Bell polynomials associated with $Y$ \par
Moving forward, our future research will focus on continuing to investigate degenerate versions and probabilistic extensions of various other special polynomials and numbers.

\end{document}